\documentclass[12pt,a4paper,reqno]{amsart}

\usepackage[latin1]{inputenc}     
\usepackage[T1]{fontenc}          
\usepackage[english]{babel}        

\usepackage{amsmath,amsfonts,amssymb,amsthm} 
\usepackage{mathtools}            
\usepackage{mathrsfs}
\usepackage{physics}

\usepackage{url}
\usepackage{cite}
\usepackage{enumitem}

\usepackage{microtype}            
\usepackage{enumitem}             
\usepackage{graphicx}             
\usepackage{booktabs}             
\usepackage{geometry}             
\usepackage[colorlinks=true, linkcolor=blue, citecolor=blue, urlcolor=blue, pdfstartview=FitH]{hyperref}

\newtheorem{theorem}{Theorem}

\newtheorem{lemma}{Lemma}

\title[Exponential stabilization of an EBB under distributed time delay]{Exponential stabilization of an Euler--Bernoulli beam via internal feedback under distributed time delay}

\author[B. B. J. Siriki]{Ben B. J. Siriki}
\address{Universit\'e Nangui Abrogoua, Abidjan, C\^ote d'Ivoire}
\email{\href{mailto: benbjsiriki@gmail.com}{benbjsiriki@gmail.com}}

\author[D. A. Behi]{Droh A. Behi}
\address{Universit\'e de Man, Man, C\^ote d'Ivoire.}
\email{arsene.behi@univ-man.edu.ci}

\author[N. F. Alla]{N'Guessan F. Alla}
\address{Universit\'e F\'elix  Houphou\"et--Boigny,
	Abidjan, C\^ote d'Ivoire.}
\email{allafabrice.nguessan@ufhb.edu.ci}

\author[A. Coulibaly]{Adama Coulibaly}
\address{Universit\'e F\'elix  Houphou\"et--Boigny,
	Abidjan, C\^ote d'Ivoire} 
\email{couliba@yahoo.fr}
\allowdisplaybreaks[4]


\begin{document}
	
	\begin{abstract}
		In this paper, we investigate the internal stabilization of a closed-loop system which consists of an Euler--Bernoulli beam subjected to a structural damping mechanism with a distributed delay into the domain. Under appropriate assumptions on the weights of the delay term and internal damping, we first prove the well-posedness of the system by a semigroup approach. Then, by constructing an appropriate Lyapunov functional and utilizing the multiplier method, we establish the exponential decay of the system energy.
	\end{abstract}
	
	\maketitle
	

	\section*{Introduction}
	
	This paper focuses on the vibration control of an Euler--Bernoulli beam model whose dynamics are governed by the initial-boundary value problem
	\begin{equation} \label{pde-system}
		\left\{ 
		\begin{aligned}
			& m y_{tt} + EI y_{xxxx} + \mathcal{U}_I = 0, \quad x \in (0,\ell), \, t>0 \\
			& y(0,t) = y_x(0,t) = y_{xx}(\ell,t) = y_{xxx}(\ell,t) = 0, \quad t>0, \\
			& y(x,0) = y_0(x), \, y_t(x,0) = y_1(x), \quad x \in [0, \ell], \\
			& y_t(x,-t) = h_0(x,-t), \quad \, x \in (0,\ell), \, t \in (0, \tau_2),
		\end{aligned} 
		\right.
	\end{equation}
	where $y(x,t)$ represents the deflection of the beam of length $\ell$ at position $x \in (0,\ell)$ and time $t > 0$. The positive physical parameters $m$ and $EI$ denote, respectively, the mass per unit length and the flexural rigidity of the beam. Furthermore, $\mathcal{U}_I(x,t)$ stands for the interior control, and the initial data $(y_0, y_1, h_0)$ belong to a suitable space.
	
	In control systems, time lags between the application of an input and the onset of the system's response are frequently observed--a phenomenon known as time delay. While very common in most dynamic processes and physical systems \cite{eurich2005, hale1993, niculescu2001} and occasionally contributing to system stability \cite{fanou2026, abdallah1993, suh1980}, delays are most often a source of instability and performance degradation \cite{datko1985, datko1991, datko1993}. Consequently, the stabilization of time-delay systems constitutes a major challenge in control theory. \newline
	Within this framework, systems involving distributed delays have recently attracted considerable attention. These models capture dynamics affected by a delay spread over a finite or infinite time interval. To counteract their destabilizing effects, various control strategies have been developed in the literature, of which we mention a few notable contributions below. 
	First, Nicaise et al. employed in \cite{nicaise2008} an observability estimate approach to establish the exponential stability of a multidimensional wave equation on a bounded domain with a smooth boundary, featuring a distributed delay on the boundary under the condition 
	\begin{equation} \label{intro:ineq1}
		\mu_0 > \int_{\tau_1}^{\tau_2} \mu(s) \, ds,
	\end{equation}
	or in the interior of the domain when 
	\begin{equation} \label{intro:ineq2}
		\mu_0 > \lVert a \rVert_{L^\infty(0,1)} \int_{\tau_1}^{\tau_2} \mu(s) \, ds,
	\end{equation}
	where $\mu_0 > 0$ is a constant, $\mu : [\tau_1, \tau_2] \to \mathbb{R}$ is a non-negative $L^\infty$ function almost everywhere, and $a$ is a positive $L^\infty$ function on the domain. We refer the reader to \cite{boumasmoud2026} for a recent application of this approach. 
	Subsequently, by combining the Sz--Nagy--Foias spectral decomposition theory \cite{szokefalvi1967} with a frequency-domain method, Saba et al. demonstrated in \cite{saba2023} the exponential stability of a microbeam equation with a distributed delay term on the boundary over a one-dimensional domain under similar hypotheses \eqref{intro:ineq1} and \eqref{intro:ineq2}. For a broader overview of this method, the reader is invited to consult \cite{silga2022, chentouf2026, zhang2022}.
	Finally, in \cite{mustafa2013}, Mustapha et al. studied a viscoelastic system affected by an internal distributed delay and subject to two types of boundary conditions. Under the smallness condition on the memory kernel $k_2$,
	\begin{equation}
		\int_{\tau_1}^{\tau_2} |k_2(s)| \, ds < k_1,
	\end{equation}
	where $k_1$ denotes the thermal diffusivity, the authors proved uniform energy decay for each configuration using a Lyapunov functional and the multiplier method \cite{komornik1994}. Recent works combining these two techniques can be found in \cite{loucif2024, mpungu2024, mpungu2021, dibes2023, feng2021, alotaibi2026}.
	
	The stabilization of Euler--Bernoulli beam equation has been the subject of active research in recent decades (see \cite{feng2019, guo2009, akil2026, azzaoui2026, siriki2026, siriki2026a, smyshlyaev2009, salhi2026, rahmtati2022}), owing to its widespread applicability in various engineering fields (see \cite{du2022, chakravarthy2010, balakrishnan1986, he2014, wright2003}). However, relatively few studies have addressed systems incorporating a distributed delay in the internal feedback law.
	This gap in the literature motivates the present work, in which we propose a feedback control consisting of a linear frictional damping coupled with a localized internal distributed delay, given by
	\begin{equation}
		\mathcal{U}_I(x,t) = \gamma_0 y_{t}(x,t) + a(x) \int_{\tau_1}^{\tau_2} \gamma(s) y_t(x,t-s) , ds, \quad x \in (0,\ell), \, t > 0,
	\end{equation}
	where $0 \leqslant \tau_1 < \tau_2$ represent the delay bounds, $\gamma_0 > 0$ is the internal damping coefficient, and the function $\gamma \in L^\infty(\tau_1, \tau_2)$ acts as the memory kernel. Furthermore, the function $a \in L^\infty(0,\ell)$ represents the spatial localization of the actuator, satisfying
	\begin{equation} \label{damp-cdt}
		\left\{
		\begin{aligned}
			& a(x) \geqslant 0 \quad \text{a.e. } x \in (0, \ell), \\
			& a(x) > a_0 > 0 \quad \text{a.e. } x \in (\omega_1, \omega_2),
		\end{aligned}
		\right.
	\end{equation}
	for any non-empty open neighbourhood $(\omega_1, \omega_2)$ of $\ell$ contained in $[0, \ell]$ .
	
	Our contribution is twofold: we establish the well-posedness of the closed-loop system \eqref{pde-system} within an appropriate energy space and prove the exponential decay of its solutions, provided that the damping coefficient and the memory kernel satisfy suitable smallness conditions. Intuitively, this demonstrates that the frictional damping effectively compensates for the adverse effects of the internal distributed delay $a(x) \int_{\tau_1}^{\tau_2} \gamma(s) y_t(x,t-s) \, ds$ to achieve exponential stabilization. 
	To the best of our knowledge, this paper provides the first stability results for an Euler--Bernoulli beam with an internal distributed delay. Despite employing a technical strategy distinct from the wave equation framework (see \cite{nicaise2008}), we establish exponential stabilization under a comparable gain condition.
	
	The rest of this paper is organized as follows. 
	Section \ref{sec1} introduces the augmented formulation of system \eqref{pde-system}, investigates its energy dissipative properties, and provide other tools required for the statement of our main results. 
	In Section \ref{sec2}, we establish the well-posedness of the system under a suitable weight condition by casting its augmented form as an abstract Cauchy problem and applying the Lumer--Phillips theorem within a semigroup framework \cite{pazy1983}. 
	Section \ref{sec3} is then dedicated to proving the exponential stability of the system by combining multiplier techniques \cite{komornik1994} with a suitable Lyapunov functional.
	
	\section{Preliminaries} \label{sec1}
	
	This section collects the preliminary material required to establish our main results. 
	Inspired by the approach developed by Nicaise and Pignotti \cite{nicaise2008}, we first define an auxiliary variable to handle the distributed delay term. Specifically, we define the new variable as:
	\begin{equation}
		u(x,\varrho,s,t) := y_t(x,t-s\varrho), \quad x \in (0,\ell), \, \varrho \in (0,1), \, s \in (\tau_1,\tau_2), \, t>0.
	\end{equation}
	Then, $u$ satisfies:
	\begin{equation}
		\begin{split}
			& s u_t + u_\varrho = 0 \quad \text{ in } (0,\ell) \times (0,1) \times (\tau_1,\tau_2) \times (0,+\infty),
			\\
			& u(x,0,s,t) = y_t(x,t), \quad  x \in (0,\ell), \, s \in (\tau_1,\tau_2), \, t>0.
		\end{split}
	\end{equation}
	Therefore, the delayed system \eqref{pde-system} is equivalent to
	\begin{gather} \label{pde-aug-system}
		\left\{
		\begin{aligned}
			& m y_{tt} + EI y_{xxxx} + \gamma_{0} u(x,0,s,t) 
			\\
			& \qquad + \int_{\tau_1}^{\tau_2} a(x) \gamma(s) u(x,1,s,t) \, ds = 0, \quad x \in (0,\ell), \, s \in (\tau_1,\tau_2), \, t>0 \\
			& s u_t + u_\varrho = 0 \quad \text{ in } (0,\ell) \times (0,1) \times (\tau_1,\tau_2) \times (0,+\infty), \\
			& y(0,t) = y_x(0,t) = y_{xx}(\ell,t) = y_{xxx}(\ell,t) = 0, \quad t>0, \\
			& u(x,0,s,t) = y_t(x,t), \quad x \in (0,\ell), \, s \in (\tau_1,\tau_2), \, t>0 \\
			& y(x,0) = y_0(x), \, y_t(x,0) = y_1(x), \quad x \in [0, \ell], \\
			& u(x,\varrho,s,0) = h_0(x,-s\varrho), \quad x \in (0, \ell), \, \varrho \in (0,1), \, s \in (0, \tau_2).
		\end{aligned} 
		\right.
	\end{gather}
	Throughout this paper, we assume that the delay and damping coefficients satisfy the condition
	\begin{align} \label{cdt:intern-coefs}
		\gamma_{0} > \lVert a \rVert_{L^\infty(0,\ell)} \int_{\tau_1}^{\tau_2} |\gamma(s)| \, ds.
	\end{align} 
	Furthermore, let us introduce the functional space $H^2_e(0,\ell)$ defined by
	\begin{equation}
		H^2_e(0,\ell) 
		:= \left\{
		y \in H^2(0,\ell) \colon y(0) = y'(0) = 0 
		\right\}.
	\end{equation}
	The estimates required for the well-posedness of system \eqref{pde-system} rely on the following norm-equivalence property.
	\begin{lemma} \label{l1}
		The map $|\cdot|_2: \, H^2_e(0,\ell) \to \mathbb{R}$ defined by
		\begin{equation} \label{norm:def}
			|w|_2 := \left( \int_0^\ell \big| w''(x) \big|^2 \, dx +   \int_0^\ell |w(x)|^2 \, dx \right)^{\frac{1}{2}}, \; \forall w \in H^2_e(0,\ell).
		\end{equation}
		defines a norm on $H^2_e(0,\ell)$ equivalent to the standard $H^2(0,\ell)$--norm. That is, there exist two positive constants $C_1, \, C_2$ such that
		\begin{gather} \label{norm:equiv}
			C_1 \lVert w \rVert_{H^2(0,\ell)} \leqslant |w|_2 \leqslant C_2 \lVert w \rVert_{H^2(0,\ell)}, \quad \forall w \in H^2_e(0,\ell).
		\end{gather}
	\end{lemma}
	
	\begin{proof}
		The verification that $| \cdot |_2$ satisfies the norm axioms on $H^2_e(0,\ell)$ is standard and omitted here. 
		Next, we prove the equivalence property. Let $w \in H^2_e(0,\ell)$. By definition, we immediately have $\lvert w \rvert_2 \leqslant \lVert w \rVert_{H^2(0,\ell)}$. 
		Furthermore, using the triangle and Cauchy--Schwarz inequalities, we get
		\begin{gather} 
			\begin{aligned}[b]
				\int_0^\ell |w'(x)|^2 \, dx 
				& \leqslant 
				\int_0^\ell x \int_0^x  |w''(\sigma)|^2 \, d\sigma \, dx
				\\
				& \leqslant \frac{\ell^2}{2} \lVert w'' \rVert^2_{L^2(0,\ell)}.
			\end{aligned}
		\end{gather}
		Consequently, we arrive at
		\begin{gather}
			\lVert w \rVert_{H^2(0,\ell)}^2 
			\leqslant \left( 1 + \frac{\ell^2}{2} \right) |w|_2^2,
		\end{gather}
		which completes the proof.
	\end{proof}
	To proceed, we define the total energy $E(t)$ of system \eqref{pde-system} by 
	\begin{align} \label{system-energy}
		\begin{split}
			E(t) 
			& := \frac{m}{2} \int_{0}^{\ell} \, y_{t}^2(x,t) \, dx + \frac{EI}{2} \int_{0}^{\ell} \, y_{xx}^2(x,t) \, dx 
			\\
			& \quad + \frac{1}{2} \int_{0}^{\ell} a(x) \int_{\tau_1}^{\tau_2} \, s |\gamma(s)| \int_{0}^{1} y_t^2(x,t-s\varrho) \, d\varrho \, ds \, dx.
		\end{split}
	\end{align}
	The following lemma shows that system \eqref{pde-system} is dissipative.
	\begin{lemma} \label{l2}
		Assume that \eqref{cdt:intern-coefs} holds. Then, for any regular solution, the energy $E(t)$ is non-increasing, namely:
		\begin{align} \label{ineq:tderiv-energy}
			\frac{dE(t)}{dt} 
			\leqslant 
			-\left( \gamma_0 - \lVert a \rVert_{L^\infty(0,\ell)} \int_{\tau_1}^{\tau_2} |\gamma(s)| \, ds \right) \int_{0}^{\ell} y_t^2 \, dx.
		\end{align}
	\end{lemma}
	
	\begin{proof}
		First, multiplying $\eqref{pde-aug-system}_1$ by $y_t$, integrating over $(0,\ell)$ by parts with respect to $x$, and taking the boundary conditions into account yields:
		\begin{equation} \label{proof-l2-4}
			\begin{split}
				& \frac{d}{dt} \left(
				\frac{m}{2} \int_{0}^{\ell} \, y_{t}^2 \, dx 
				+ \frac{EI}{2} \int_{0}^{\ell} \, y_{xx}^2 \, dx \right)
				\\
				& \qquad 
				= - \gamma_0 \int_{0}^{\ell} y_t^2 \, dx 
				- \frac{\xi}{2} \int_{0}^{\ell} a(x) \int_{\tau_1}^{\tau_2} \gamma(s) u(x,1,s,t) \, y_t \, ds \, dx.
			\end{split}
		\end{equation}
		Next, multiplying $\eqref{pde-aug-system}_2$ by $a(x)\vert{}\gamma(s)\vert{} u_\varrho$ and integrating by parts with respect to $\varrho$ gives:
		\begin{equation} \label{proof-l2-7}
			\begin{aligned}
				& \frac{d}{dt} \left( \int_{0}^{\ell} a(x) \int_{\tau_1}^{\tau_2} s |\gamma(s)| \, \int_{0}^{1} \, u^2 \, d\varrho \, ds \, dx \right)
				\\
				& \quad =
				\int_{0}^{\ell} a(x) \int_{\tau_1}^{\tau_2} |\gamma(s)| \left( y_t^2 - u(x,1,s,t)^2 \right) \, ds \, dx.
			\end{aligned}	
		\end{equation}
		Then, summing expressions \eqref{proof-l2-4} and \eqref{proof-l2-7} leads to:
		\begin{equation} \label{proof-l2-8}
			\begin{aligned}[b]
				\frac{dE(t)}{dt} 
				& \quad
				= \int_{0}^{\ell} \left( -\gamma_0 + \frac{1}{2} a(x) \int_{\tau_1}^{\tau_2} |\gamma(s)| \, ds \right) y_t^2 \, dx  
				\\
				& \qquad 
				- \frac{1}{2} \int_{0}^{\ell} a(x) \int_{\tau_1}^{\tau_2} |\gamma(s)| \, u(x,1,s,t)^2 \, ds \, dx
				\\
				& \qquad 
				- \int_{0}^{\ell} a(x) \int_{\tau_1}^{\tau_2} \gamma(s) u(x,1,s,t) \, y_t \, ds \, dx.
			\end{aligned}
		\end{equation}
		Applying Young's inequality to the integrand $-\gamma(s) u(x,1,s,t) y_t$, we obtain:
		\begin{align} \label{proof-l2-9}
			\left| \gamma(s) u(x,1,s,t) y_t \right|
			\leqslant \frac{|\gamma(s)|}{2} \Big( u(x,1,s,t)^2 + y_t^2 \Big).
		\end{align}
		Finally, inequality \eqref{ineq:tderiv-energy} follows directly by substituting \eqref{proof-l2-9} into \eqref{proof-l2-8} and utilizing condition \eqref{cdt:intern-coefs}. 
	\end{proof}
	This section concludes with two key estimates stated in the following lemma, which will play a crucial role in the proof of the stability result.
	\begin{lemma} \label{l3}
		Let $y \in H^2_e(0,\ell)$. Then the following estimates hold:
		\begin{gather}
			\int_0^\ell y^2(x) \, dx 
			\leqslant
			\frac{\ell^4}{8} \int_0^\ell \, \big( y''(x) \big)^2 \, dx 
			\label{ineq1-Poincare}
			\\
			\int_0^\ell \, a(x) \, y^2(x) \, dx 
			\leqslant 
			\frac{\ell^3}{2} \int_0^\ell \, a(x) \, dx \, \big\lVert y'' \big\rVert_{L^2(0,\ell)}^2. 
			\label{ineq2-Poincare}
		\end{gather}
	\end{lemma}
	
	\begin{proof}
		We have for all $x \in (0,\ell)$:
		\begin{gather} \label{proof-l3-1}
			y^2(x) 
			= \left( \int_0^x \, y'(r) \, dr \right)^2 
			\leqslant 
			x \int_0^x \, y'(r)^2 \, dr.
		\end{gather}
		Analogously, we have for all $x \in (0,\ell)$:
		\begin{gather} \label{proof-l3-2}
			y'^2(x) 
			\leqslant 
			x \big\lVert y'' \big\rVert_{L^2(0,\ell)}^2.
		\end{gather}
		Then, we deduce from inequalities \eqref{proof-l3-1} and \eqref{proof-l3-2} that
		\begin{gather*}
			\begin{aligned}
				\int_0^\ell y^2(x) \, dx 
				& \leqslant 
				\int_0^\ell x \int_0^x \, r \, dr \, dx \; \big\lVert y'' \big\rVert_{L^2(0,\ell)}^2 
				=
				\frac{\ell^4}{8} \, \big\lVert y'' \big\rVert_{L^2(0,\ell)}^2
			\end{aligned}
			\text{and }
			\\
			\begin{aligned}[b]
				\int_0^\ell a(x) y^2(x) \, dx 
				& \leqslant 
				\int_0^\ell x a(x) \int_0^x \, r \, dr \, dx \, \big\lVert y'' \big\rVert_{L^2(0,\ell)}^2 
				=
				\frac{\ell^3}{2} \int_0^\ell  a(x) \, dx \, \big\lVert y'' \big\rVert_{L^2(0,\ell)}^2.
				\\
			\end{aligned}
		\end{gather*}
	\end{proof}
	
	\section{Existence and uniqueness of solution} \label{sec2}
	
	This section addresses the well-posedness of system \eqref{pde-system}. To this end, the Lumer--Phillips theorem is applied to an abstract Cauchy problem associated with the augmented system \eqref{pde-aug-system}. 
	Without loss de generality, we assume throughout this section that all functional spaces are complex-valued.
	%
	
	Consider the following Hilbert space 
	\begin{align}
		\mathscr{H} := H^2_e(0,\ell) \times L^2(0, \ell) \times L^2\big((0,\ell) \times (\tau_1,\tau_2) \times (0,1)\big),
	\end{align}
	equipped with the following inner product:
	\begin{align}
		\begin{split}
			\langle \Psi_1, \Psi_2 \rangle 
			&:= 
			m \int_{0}^{\ell} v_1 \overline{v_2} \, dx + EI \int_{0}^{\ell} \, y_1'' \overline{y_2''} dx 
			\\
			& \quad + \int_{0}^{\ell} a(x) \int_{\tau_1}^{\tau_2} s |\gamma(s)| \int_{0}^{1} u_1 \overline{u_2} \, d\varrho \, ds \, dx,
		\end{split}
	\end{align}
	for all $\Psi_i = (y_i, v_i, u_i) \in \mathscr{H}, \, i = 1, 2$. Using the semigroup approach \cite{pazy1983}, the system \eqref{pde-aug-system} can be formulated as an abstract Cauchy problem:
	\begin{gather}
		\left\{
		\begin{aligned}
			& \frac{d}{dt} \Psi(t) = \mathcal{A} \Psi(t), 
			\\
			& \Psi(0) = \Psi_0 = \big( y_0, \, y_1, \, h_0 \big),
		\end{aligned}
		\right.
	\end{gather}
	where $\Psi(t) := (y(\cdot,t), v(\cdot,t), u(\cdot,\varrho,s,t))$ and the operator $\mathcal{A}$ with its domain $\mathscr{D}(\mathcal{A})$ are defined by:
	\begin{gather}
		\begin{gathered}
			\mathscr{D}\big(\mathcal{A}\big) := 
			\left\{
			\begin{pmatrix}
				y \\ v \\ u
			\end{pmatrix} \in \mathscr{H}: \, 
			\begin{matrix}
				y \in H^4(0,\ell), \, v \in H^2_e(0,\ell), 
				\\
				u \in L^2\big( (0,\ell) \times (\tau_1,\tau_2); H^1(0,1) \big) 
				\\
				u(\cdot,0,s) = v, \, y(0) = y'(0) = y''(\ell) = y^{(3)}(\ell) = 0
			\end{matrix}
			\right\},
		\end{gathered}
		\\
		\mathcal{A} \Psi :=
		\begin{pmatrix}
			v \\ 
			-\frac{EI}{m} y^{(4)} - \frac{1}{m} \Big( \gamma_0 v + a(\cdot) \int_{\tau_1}^{\tau_2} \gamma(s) u(\cdot,1,s) \, ds \Big) \\
			-s^{-1} u_\varrho
		\end{pmatrix}.
	\end{gather}
	
	We now state and prove the well-posedness of system \eqref{pde-system}.
	\begin{theorem} \label{t1}
		Assume that \eqref{cdt:intern-coefs} holds. For any initial value $\Psi_0 \in \mathscr{H}$, system \eqref{pde-system} has a unique weak solution $\Psi \in C\left(\mathbb{R}_+, \mathscr{H}\right)$. 
		Moreover, if $\Psi_0 \in \mathscr{D}(\mathcal{A})$, then $\Psi \in C\left(\mathbb{R}_+, \mathscr{D}(\mathcal{A})\right) \cap C^1\left(\mathbb{R}_+, \mathscr{H}\right)$.
	\end{theorem}
	
	\begin{proof}
		By Lumer--Phillips theorem, it suffices to show that the linear operator $\mathcal{A}$ is $m$-dissipative. First, to establish the dissipativeness of $\mathcal{A}$, let $\Psi \in \mathscr{D}(\mathcal{A})$. By definition,
		\begin{equation} \label{proof-t1-1}
			\begin{aligned}
				\langle \mathcal{A}\Psi, \Psi \rangle 
				&=
				- \gamma_0 \int_{0}^{\ell} |v(x)|^2 \, dx - EI \int_{0}^{\ell} y^{(4)} \overline{v} \, dx + EI \int_{0}^{\ell} v'' \overline{y''} \, dx 
				\\
				& \quad
				- \int_{0}^{\ell} a(x) \int_{\tau_1}^{\tau_2} \gamma(s) u(x,1,s) \overline{v} \, ds \, dx 
				- \int_{0}^{\ell} a(x) \int_{\tau_1}^{\tau_2} |\gamma(s)| \int_{0}^{1} u_\varrho \overline{u} \, d\varrho \, ds \, dx.
			\end{aligned}
		\end{equation}
		Integrating by parts the second and fifth terms of the above expression and taking the boundary conditions $\eqref{pde-aug-system}_3$ and $\eqref{pde-aug-system}_4$ into account, we obtain:
		\begin{equation} \label{proof-t1-4}
			\begin{aligned}[b]
				\langle \mathcal{A}\Psi, \Psi \rangle 
				&=
				- \gamma_0 \int_{0}^{\ell} |v(x)|^2 \, dx 
				- \int_{0}^{\ell} a(x) \int_{\tau_1}^{\tau_2} \gamma(s) u(x,1,s) \, ds \, \overline{v} \, dx  
				\\
				& \quad 
				- \int_{0}^{\ell} a(x) \int_{\tau_1}^{\tau_2} |\gamma(s)| \int_{0}^{1} u \overline{u_\varrho} \, d\varrho \, ds \, dx 
				+ 2 i EI \Im \left( \int_{0}^{\ell} v'' \overline{y''} \, dx \right),
			\end{aligned}
		\end{equation}
		which implies that
		\begin{equation}
			\begin{aligned} \label{proof-t1-4a}
				\Re \langle \mathcal{A} \Psi, \Psi \rangle 
				&=
				\int_{0}^{\ell} \left( - \gamma_0 + \frac{a(x)}{2} \int_{\tau_1}^{\tau_2} |\gamma(s)| \, ds \right) |v(x)|^2 \, dx 
				\\
				& \quad 
				- \frac{a(x)}{2}\int_{0}^{\ell} \int_{\tau_1}^{\tau_2} |\gamma(s)| |u(x,1,s)|^2 \, ds \, dx
				\\
				& \quad 
				- \int_{0}^{\ell} a(x) \int_{\tau_1}^{\tau_2} \gamma(s) \Re \left( u(x,1,s) \overline{v} \right) \, ds \, dx.
			\end{aligned}
		\end{equation}
		By applying Young's inequality on the integrand of the last integral term of the precedent identity, we get:
		\begin{align} \label{proof-t1-5}
			-\gamma(s) \Re \left( u(\cdot,1,s) \overline{v} \right)
			\leqslant
			\frac{|\gamma(s)|}{2} \left( |v|^2 + |u(\cdot,1,s)|^2 \right).
		\end{align}
		Then, substituting this estimate into \eqref{proof-t1-4} and utilizing conditions \eqref{cdt:intern-coefs} and $\eqref{damp-cdt}_1$ leads to:
		\begin{equation}
			\begin{aligned}[b]
				\Re \langle \mathcal{A} \Psi, \Psi \rangle 
				&\leqslant
				\int_{0}^{\ell} \left( - \gamma_0 + a(x) \int_{\tau_1}^{\tau_2} |\gamma(s)| \, ds \right) |v(x)|^2 \, dx
				\\
				&\leqslant
				\left( -\gamma_0 + \lVert a \rVert_{L^\infty(0,\ell)} \int_{\tau_1}^{\tau_2} |\gamma(s)| \, ds \right) \int_{0}^{\ell} |v(x)|^2 \, dx \leqslant 0.
			\end{aligned}
		\end{equation}
		Next, for the maximality of $\mathcal{A}$, we demonstrate that $\lambda I - \mathcal{A}$ is surjective for some $\lambda > 0$. Given $F = (f,g,h) \in \mathscr{H}$, this amounts to showing that the equation $(\lambda I - \mathcal{A})\Psi = F$ admits a unique solution $\Psi \in \mathscr{D}(\mathcal{A})$. Equivalently, the resolvent equation reads:
		\begin{gather} \label{proof-t1-6}
			\left\{\begin{aligned}
				& \lambda y - v = f \\
				& \lambda v + \frac{EI}{m} y^{(4)} + \frac{1}{m} \left( \gamma_0 v + a(\cdot) \int_{\tau_1}^{\tau_2} \gamma(s) u(\cdot,1,s) \, ds \right) = g \\
				& \lambda u + s^{-1} u_\varrho = h.
			\end{aligned} 
			\right.
		\end{gather}
		Solving equations $\eqref{proof-t1-6}_1$ and $\eqref{proof-t1-6}_3$ for $\varrho$ and using the fact that $u(\cdot,0,s) = v$, we obtain:
		\begin{align}
			& v = \lambda y - f, \label{proof-t1-7}
			\\
			& u(\cdot,\varrho,s) = e^{-\lambda s \varrho} \left( v + s \int_{0}^{\varrho} h(\cdot,r,s) e^{s \lambda r} \, dr \right). \label{proof-t1-8}
		\end{align}
		Thus, the solvability of system \eqref{proof-t1-6} reduces to finding $y$. Upon incorporating these expressions into $\eqref{proof-t1-6}_2$, we obtain the following system for $y$: 
		\begin{gather} \label{proof-t1-9}
			\left\{
			\begin{aligned}
				& EI y^{(4)} + \lambda \, \Xi_{1}(\lambda,\cdot) y = m g + \Xi_{1}(\lambda,\cdot) f + \Xi_{2}h(\lambda,\cdot) \\
				& y(0) = y'(0) = y''(\ell) = y^{(3)}(\ell) = 0,
			\end{aligned}
			\right. 
		\end{gather}
		where the functions $\Xi_{i}, \, i = 1,2$ are defined for all $x \in (0,\ell)$ by:
		\begin{gather} 
			\Xi_{1}(\lambda,x) := m \lambda + \gamma_{0} + a(x) \int_{\tau_1}^{\tau_2} \gamma(s) e^{-\lambda s} \, ds,
			\\
			\Xi_{2}h(\lambda,x) := - a(x) \int_{\tau_1}^{\tau_2} s \gamma(s) e^{-\lambda s} \int_{0}^{1} h(x,r,s) e^{\lambda s r} \, dr \, ds.
		\end{gather}
		By virtue of condition \eqref{damp-cdt}, the non-negativity of $\Xi_1$ is readily established. \newline
		To derive the variational formulation, we multiply equation $\eqref{proof-t1-9}_1$ by $\overline{\varphi}$, integrate by parts over $(0,\ell)$, and apply boundary condition $\eqref{pde-system}_2$. This yields the weak problem, which requires finding $y \in H^2_e(0,\ell)$ such that: 
		\begin{align} \label{proof-t1-11}
			\mathcal{Q}(y,\varphi) = \mathcal{L}(\varphi),
		\end{align}
		for all $\varphi \in H^2_e(0,\ell)$, where the sesquilinear form $\mathcal{Q}$ and the antilinear form $\mathcal{L}$ are defined by:
		\begin{gather}
			\mathcal{Q}(y,\varphi) 
			= EI \int_0^\ell y'' \overline{\varphi''} \, dx + \lambda \int_0^\ell \Xi_{1}(\lambda,x) \, y  \overline{\varphi} \, dx
			\\ 
			\mathcal{L}(\varphi)
			= \int_0^\ell  \left( m g + \Xi_{1}(\lambda,x) f + \Xi_{2}h(\lambda,x) \right) \overline{\varphi} \, dx. 
		\end{gather}
		Direct calculations based on the Cauchy--Schwarz and triangle inequalities, Sobolev embeddings, and conditions \eqref{cdt:intern-coefs}--\eqref{damp-cdt} establish the continuity of $\mathcal{Q}$ and $\mathcal{L}$.
		Furthermore, invoking the norm-equivalence property (see Lemma \ref{l1}), we deduce that $\mathcal{Q}$ is coercive:
		\begin{align}
			\Re \mathcal{Q}(y,y) 
			& = EI \int_0^\ell \big| y'' \big|^2 \, dx + \lambda \int_0^\ell \Xi_{1}(\lambda,x) \, |y|^2 \, dx \nonumber
			\\
			\Re \mathcal{Q}(y,y) 
			& \geqslant 
			\min \left\{ EI; m \lambda^2 \right\} C_1 \lVert y \rVert^2_{H^2(0,\ell)}.
		\end{align}
		According to Lax-Milgram theorem, the system \eqref{proof-t1-9} admits a unique weak solution $y \in H^2_e(0,\ell)$. Consequently, $v \in H^2_e(0,\ell)$ and $u \in L^2\big((0,\ell) \times (\tau_1,\tau_2) \times (0,1)\big)$.
		
		Conversely, let $y \in H^2_e(0,\ell)$ be a weak solution. Integrating by parts, equation \eqref{proof-t1-11} takes the form:
		\begin{gather} \label{proof-t1-17}
			\begin{aligned}
				& \int_0^\ell \left( EI y^{(4)} + \lambda \, \Xi_{1} (\lambda,x) \, y  \right) \overline{\varphi} \, dx 
				+ EI y''(\ell) \overline{\varphi'(\ell)} - EI y^{(3)} (\ell) \overline{\varphi(\ell)} 
				\\
				& \qquad \quad 
				= \int_0^\ell  \Big( m g + \Xi_{1} (\lambda,x) f + \Xi_{2}h (\lambda,x) \Big) \overline{\varphi} \, dx.
			\end{aligned}
		\end{gather}
		First, letting $\varphi \in \mathscr{D}(0,\ell)$, the preceding identity becomes:
		\begin{equation}
			\begin{gathered}
				\int_0^\ell \left( EI y^{(4)} + \lambda \, \Xi_{1} (\lambda,x) \, y  \right) \overline{\varphi} \, dx 
				= \int_0^\ell  \Big( m g + \Xi_{1} (\lambda,x) f + \Xi_{2}h (\lambda,x) \Big) \overline{\varphi} \, dx.
			\end{gathered}
		\end{equation}
		This implies that $EI y^{(4)} + \lambda \Xi_{1} (\lambda,\cdot) \, y = m g + \Xi_{1} (\lambda,\cdot) f + \Xi_{2}h (\lambda,\cdot)$ almost everywhere on $(0,\ell)$. Hence, $y \in H^4(0,\ell)$ and satisfies equation $\eqref{proof-t1-9}_1$ by virtue of the density of $\mathscr{D}(0,\ell)$ in $L^2(0,\ell)$. \newline
		Alternatively, let $\varphi \in H^2_e(0,\ell)$. In view of the preceding results, equation \eqref{proof-t1-17} yields:
		\begin{equation}
			EI y''(\ell) \overline{\varphi'(\ell)} - EI y^{(3)} (\ell) \overline{\varphi(\ell)} = 0.
		\end{equation}
		Due to the arbitrariness of $\varphi$, the solution $y$ satisfies $\eqref{proof-t1-9}_2$. 
		Thus, $y$ is a solution to \eqref{proof-t1-9}, which establishes the maximality of the operator $\mathcal{A}$.
	\end{proof}
	
	\section{Uniform decay result} \label{sec3}
	This section is devoted to the exponential stability of system \eqref{pde-system}. Our approach relies on the construction of a suitable Lyapunov functional combined with multiplier techniques.
	
	\subsection{Technical lemmas}
	To establish the uniform decay rate of the energy $E(t)$ associated with system \eqref{pde-system}, we establish several technical lemmas.
	We introduce the functional $F(t)$ defined by:
	\begin{gather} \label{aux-funct}
		F(t) := \zeta_1 F_1(t) + \zeta_2 F_2(t), \quad t \geqslant 0,
	\end{gather}
	where the constants $\zeta_1, \, \zeta_2>0$ and $F_i, \, i=1,2$ are given by:
	\begin{gather}
		F_1(t) := m \int_{0}^{\ell} y y_t \, dx
		\\
		F_2(t) := \int_{0}^{\ell} a(x) \int_{\tau_1}^{\tau_2} \, s |\gamma(s)| \int_{0}^{1} e^{-s\varrho} y_t^2(x,t-s\varrho) \, d\varrho \, ds \, dx.
	\end{gather}
	
	\begin{lemma} \label{l4}
		Let $y$ be the solution of system \eqref{pde-system}. Then, for every $\varepsilon_1 > 0$, the following estimate holds:
		\begin{gather} \label{estim:tderiv-f1}
			\begin{aligned}
				F_1'(t) 
				& \leqslant 
				\left( \frac{m}{2} + \frac{\varepsilon_1 \gamma_0}{2} \right) \int_{0}^{\ell} y_t^2 \, dx
				+ \frac{\varepsilon_1}{2} \int_{0}^{\ell} a(x) \int_{\tau_1}^{\tau_2} \big| \gamma(s) \big| \, y_t^2(x,t-s) \, ds \, dx
				\\
				& \quad 
				+ \left[ -EI +  \frac{\ell^4}{4\varepsilon_1} \left( \frac{\gamma_0}{4} + \lVert a \rVert_{L^\infty(0,\ell)} \int_{\tau_1}^{\tau_2} \big| \gamma(s) \big| \, ds \right) \right] \int_0^\ell y_{xx}^2 \, dx.
			\end{aligned}
		\end{gather}
	\end{lemma}
	
	\begin{proof}
		Differentiating with respect to time and utilizing equation $\eqref{pde-system}_1$, we obtain:
		\begin{gather} \label{proof-l4-1}
			\begin{aligned}
				F_1'(t) 
				& =
				m \int_{0}^{\ell} y_t^2 \, dx 
				- EI \int_{0}^{\ell} y y_{xxxx} \, dx 
				\\
				& \quad 
				- \int_{0}^{\ell} y \left( \gamma_{0} y_{t} + \int_{\tau_1}^{\tau_2} a(x) \gamma(s) y_t(x,t-s) \, ds \right) \, dx.
			\end{aligned}
		\end{gather}
		Integrating the second term on the right-hand side by parts and invoking boundary condition $\eqref{pde-system}_2$, relation \eqref{proof-l4-1} becomes:
		\begin{equation} \label{proof-l4-3}
			\begin{split}
				F_1'(t) 
				& = 
				\frac{m}{2} \int_{0}^{\ell} y_t^2 \, dx
				- EI \int_{0}^{\ell} y_{xx}^2 \, dx 
				\\
				& \quad 
				- \int_{0}^{\ell} \left( \gamma_{0} y_{t} + \int_{\tau_1}^{\tau_2} a(x) \gamma(s) y_t(x,t-s) \, ds \right) y \, dx.
			\end{split}
		\end{equation}
		Next, we estimate the last integral term in the preceding equality. Applying the triangle inequality, we obtain:
		\begin{gather} \label{proof-l4-4}
			\begin{aligned}
				& \left| \int_{0}^{\ell} \left( \gamma_{0} y_{t} + \int_{\tau_1}^{\tau_2} a(x) \gamma(s) y_t(x,t-s) \, ds \right) y \, dx \right| 
				\\
				& \quad \leqslant 
				\gamma_0 \int_{0}^{\ell} \big| y y_{t} \big| \, dx 
				+ \int_{\tau_1}^{\tau_2} \int_{0}^{\ell} a(x)  \big| \gamma(s) \big| \big| y y_t(x,t-s) \big| \, ds \, dx.
			\end{aligned}
		\end{gather}
		Now, applying Young's inequality to each term on the right-hand side of \eqref{proof-l4-4} and using inequalities \eqref{ineq1-Poincare} and \eqref{ineq2-Poincare} (see Lemma \ref{l3}) yields, for all $\varepsilon_1 > 0$:
		\begin{gather} \label{proof-l4-5}
			\begin{aligned}
				& \int_{0}^{\ell} \big| y y_{t} \big| \, dx 
				\leqslant
				\frac{\varepsilon_1}{2} \int_{0}^{\ell} y_t^2 \, dx + \frac{\ell^4}{16\varepsilon_1} \int_{0}^{\ell} y_{xx}^2 \, dx;
			\end{aligned}
			\\
			\label{proof-l4-5a}
			\begin{aligned}[b]
				& \int_{\tau_1}^{\tau_2} \int_{0}^{\ell} a(x)  \big| \gamma(s) \big| \big| y y_t(x,t-s) \big| \, ds \, dx
				\\
				& \quad \leqslant 
				\frac{\varepsilon_1}{2} \int_{0}^{\ell} a(x) \int_{\tau_1}^{\tau_2}  \big| \gamma(s) \big| y_t^2(x,t-s) \, ds \, dx 
				+ \frac{1}{2\varepsilon_1} \int_{\tau_1}^{\tau_2} \big| \gamma(s) \big| \, ds \, \int_{0}^{\ell} a(x) y^2 \, dx
				\\
				& \quad \leqslant 
				\frac{\varepsilon_1}{2} \int_{0}^{\ell} a(x) \int_{\tau_1}^{\tau_2}  \big| \gamma(s) \big| y_t^2(x,t-s) \, ds \, dx 
				+ \frac{\ell^3}{4\varepsilon_1} \int_{\tau_1}^{\tau_2} \big| \gamma(s) \big| \, ds \, \int_0^\ell \, a(x) \, dx \, \int_0^\ell y_{xx}^2 \, dx.
			\end{aligned}
		\end{gather}
		Inserting \eqref{proof-l4-5} and \eqref{proof-l4-5a} into \eqref{proof-l4-4}, we deduce that for any $\varepsilon_1 > 0$:
		\begin{gather} \label{proof-l4-6}
			\begin{aligned}[b]
				& \left| \int_{0}^{\ell} \left( \gamma_{0} y_{t} + \int_{\tau_1}^{\tau_2} a(x) \gamma(s) y_t(x,t-s) \, ds \right) y \, dx \right| 
				\\
				& \quad \leqslant 
				\gamma_0 \int_{0}^{\ell} \left( \frac{\varepsilon_1}{2} y_t^2 + \frac{\ell^4}{16\varepsilon_1} y_{xx}^2 \right) \, dx 
				+ \frac{\ell^3}{4\varepsilon_1} \int_{\tau_1}^{\tau_2} \big| \gamma(s) \big| \, ds \, \int_0^\ell \, a(x) \, dx \, \int_0^\ell y_{xx}^2 \, dx
				\\
				& \qquad \quad 
				+ \frac{\varepsilon_1}{2} \int_{0}^{\ell} a(x) \int_{\tau_1}^{\tau_2} \big| \gamma(s) \big| \, y_t^2(x,t-s) \, ds \, dx
				\\
				& \quad \leqslant 
				\frac{\varepsilon_1 \gamma_0}{2}  \int_{0}^{\ell} y_t^2 \, dx 
				+ \frac{\varepsilon_1}{2} \int_{0}^{\ell} a(x) \int_{\tau_1}^{\tau_2} \big| \gamma(s) \big| \, y_t^2(x,t-s) \, ds \, dx
				\\
				& \qquad 
				+ \frac{\ell^4}{4\varepsilon_1} \left( \frac{\gamma_0}{4} + \lVert a \rVert_{L^\infty(0,\ell)} \int_{\tau_1}^{\tau_2} \big| \gamma(s) \big| \, ds \right) \int_0^\ell y_{xx}^2 \, dx.
			\end{aligned}
		\end{gather}
		%
		Finally, the desired estimate follows directly by substituting \eqref{proof-l4-6} into \eqref{proof-l4-3}.
	\end{proof}
	
	\begin{lemma} \label{l5}
		Let $y$ be the solution of system \eqref{pde-system}. Then, the time derivative of $F_2$ satisfies:
		\begin{gather} \label{estim:tderiv-f2}
			\begin{aligned}[b]
				F_2'(t) 
				& \leqslant  
				\int_{\tau_1}^{\tau_2} \, |\gamma(s)| \, ds \, \int_{0}^{\ell} a(x) y_{t}^2(x,t) \, dx 
				\\
				& \quad - e^{-\tau_2} \int_{0}^{\ell} a(x) \int_{\tau_1}^{\tau_2} \, |\gamma(s)| \, y_{t}^2(x,t-s) \, ds \, dx
				\\
				& \quad - e^{-\tau_2} \int_{0}^{\ell} a(x) \int_{\tau_1}^{\tau_2} \, |\gamma(s)| \int_{0}^{1} y_{t}^2(x,t-s\varrho) \, d\varrho \, ds \, dx.
			\end{aligned}
		\end{gather}
	\end{lemma}
	
	\begin{proof}
		Combining equation $\eqref{pde-aug-system}_2$ with integration by parts, we get:
		\begin{gather*}
			\begin{aligned}[b]
				F_2'(t) 
				& = 2 \int_{0}^{\ell} a(x) \int_{\tau_1}^{\tau_2} \, s |\gamma(s)| \int_{0}^{1} e^{-s\varrho} y_t(x,t-s\varrho) y_{tt}(x,t-s\varrho) \, d\varrho \, ds \, dx
				\\ 
				& = - \int_{0}^{\ell} a(x) \int_{\tau_1}^{\tau_2} \, |\gamma(s)| \int_{0}^{1} e^{-s\varrho} \frac{\partial}{\partial s} \Big( y_{t}^2(x,t-s\varrho) \Big) \, d\varrho \, ds \, dx
				\\ 
				& = 
				\int_{\tau_1}^{\tau_2} \, |\gamma(s)| \, ds \, \int_{0}^{\ell} a(x) y_{t}^2(x,t) \, dx 
				- \int_{0}^{\ell} a(x) \int_{\tau_1}^{\tau_2} \, |\gamma(s)| e^{-s} y_{t}^2(x,t-s) \, ds \, dx
				\\
				& \quad - \int_{0}^{\ell} a(x) \int_{\tau_1}^{\tau_2} \, s|\gamma(s)| \int_{0}^{1} e^{-s\varrho} y_{t}^2(x,t-s\varrho) \, d\varrho \, ds \, dx.
			\end{aligned}
		\end{gather*}
		Hence, in view of the inequality $-e^{-s\varrho} \leqslant -e^{-s} \leqslant -e^{-\tau_2} < 0$ for all $(s,\varrho) \in [\tau_1, \tau_2] \times [0,1]$, the desired estimate follows. 
	\end{proof}
	
	We are now in a position to estimate the time derivative of the functional $F(t)$. 
	\begin{lemma} \label{l6}
		Assume that condition \eqref{cdt:intern-coefs} holds, and let $y$ be the solution to the system \eqref{pde-system}. Then, for all $t \geqslant 0$ and $\varepsilon_1>0$,
		\begin{gather} \label{estim:tderiv-aux-funct}
			\begin{aligned}[b]
				F'(t)
				& \leqslant
				\int_{0}^{\ell} \left( \left( \frac{m}{2} + \frac{\varepsilon_1 \gamma_0}{2} \right) \zeta_1 + \zeta_2 a(x) \int_{\tau_1}^{\tau_2} \, |\gamma(s)| \, ds \right) y_t^2 \, dx 
				\\ 
				& \quad 
				+ \zeta_1 \left[ -EI +  \frac{\ell^4}{4\varepsilon_1} \left( \frac{\gamma_0}{4} + \lVert a \rVert_{L^\infty(0,\ell)} \int_{\tau_1}^{\tau_2} \big| \gamma(s) \big| \, ds \right) \right] \int_0^\ell y_{xx}^2 \, dx 
				\\
				& \quad 
				- e^{-\tau_2} \zeta_2 \int_{0}^{\ell} a(x) \int_{\tau_1}^{\tau_2} \, |\gamma(s)| \int_{0}^{1} y_{t}^2(x,t-s\varrho) \, d\varrho \, ds \, dx 
				\\
				& \quad 
				+ \left( \frac{\varepsilon_1}{2} \zeta_1 - e^{-\tau_2} \zeta_2 \right) \int_{0}^{\ell} a(x) \int_{\tau_1}^{\tau_2} \big| \gamma(s) \big| \, y_t^2(x,t-s) \, ds \, dx.
			\end{aligned}
		\end{gather}
	\end{lemma}
	
	\begin{proof}
		The result is immediate by combining estimates \eqref{estim:tderiv-f1} and \eqref{estim:tderiv-f2} from Lemmas \ref{l4}--\ref{l5}.
	\end{proof}

	\subsection{Exponential stability}
	
	We define the Lyapunov functional by 
	\begin{align}
		\mathscr{L}(t) := E(t) + F(t), \quad t \geqslant 0,
	\end{align}
	where $E(t)$ and $F(t)$ are given by \eqref{system-energy} and \eqref{aux-funct}, respectively. \newline 
	The equivalence between the functional $\mathscr{L}(t)$ and the system energy $E(t)$ is established in the following lemma. 
	\begin{lemma} \label{l7}
		Let $\zeta_1 > 0$ be sufficiently small. Then there exist two positive constants $\mu_1$ and $\mu_2$ such that, for any $t \geqslant 0$,
		\begin{gather} \label{equiv-energy-lyap}
			\mu_1 E(t) \leqslant \mathscr{L}(t) \leqslant \mu_2 E(t).
		\end{gather}
	\end{lemma}

	\begin{proof}
		Using Young's inequality and \eqref{ineq1-Poincare}, we obtain, for every $\varepsilon_2 > 0$: 
		\begin{equation} 
			\begin{aligned}[b]
				|F_1(t)| 
				& \leqslant
				m \left( \frac{\varepsilon_2}{2} \int_{0}^{\ell} y^2(x,t) \, dx + \frac{1}{2\varepsilon_2} \int_{0}^{\ell} y_t^2(x,t) \, dx \right)
				\\
				& \leqslant 
				m \left( \frac{\varepsilon_2 \ell^4}{16}  \lVert y_{xx} \rVert_{L^2(0,\ell)}^2 + \frac{1}{2\varepsilon_2} \lVert y_t \rVert_{L^2(0,\ell)}^2 \right).
			\end{aligned}
		\end{equation}
		Taking $\varepsilon_2 = \frac{2\sqrt{2}}{\ell^2}$, the previous inequality simplifies to: 
		\begin{gather} \label{proof-l7-1}
			\begin{aligned}[b]
				|F_1(t)| 
				& \leqslant 
				\frac{m \ell^2}{4 \sqrt{2}} \left( \lVert y_{xx} \rVert_{L^2(0,\ell)}^2 + \lVert y_t \rVert_{L^2(0,\ell)}^2 \right) 
				\\
				& \leqslant C E(t),
			\end{aligned}
		\end{gather}
		where the positive constant $C := \frac{\ell^2}{2\sqrt{2}} \max \left\{ 1, \, \frac{m}{EI}  \right\}$. Next, we have:
		\begin{equation} \label{proof-l7-2}
			\begin{aligned}[b]
				F_2(t)
				& \leqslant
				\int_{0}^{\ell} a(x) \int_{\tau_1}^{\tau_2} \, s |\gamma(s)| \int_{0}^{1} y_t^2(x,t-s\varrho) \, d\varrho \, ds \, dx 
				\\
				& \leqslant 2 E(t).
			\end{aligned}
		\end{equation}
		Hence, combining estimates \eqref{proof-l7-1} and \eqref{proof-l7-2} yields
		\begin{gather*}
			\begin{aligned}
				\left( 1 - \zeta_1 C \right) E(t) \leqslant \mathscr{L}(t) \leqslant \left( 1 + 2 \zeta_2 + \zeta_1 C \right) E(t),
			\end{aligned}
		\end{gather*} 
		provided that $\zeta_1$ is chosen such that 
		\begin{gather} \label{proof-l7-3}
			0 < \zeta_1 < C^{-1}.
		\end{gather}
	\end{proof}
	
	We are now in a position to prove the main result of this paper.
	\begin{theorem}
		Assume that condition \eqref{cdt:intern-coefs} holds, and let $y$ be the solution to the system \eqref{pde-system}. Then, there exist two positive constants $C$ and $\omega$ such that
		\begin{align} \label{ineq:expo-decay}
			E(t) \leqslant C e^{-\omega t}, \quad \forall t \geqslant 0.
		\end{align}
	\end{theorem}

	\begin{proof}
		Using the time derivative estimates \eqref{ineq:tderiv-energy} and \eqref{estim:tderiv-aux-funct} provided by Lemmas \ref{l2} and \ref{l6}, we obtain:
		\begin{gather} \label{proof-t2-1}
			\begin{aligned}[b]
				\mathscr{L}'(t) 
				& = E'(t) + F'(t) 
				\\
				& \leqslant 
				-\left( \gamma_0 - \lVert a \rVert_{L^\infty(0,\ell)} \int_{\tau_1}^{\tau_2} |\gamma(s)| \, ds \right) \int_{0}^{\ell} y_t^2 \, dx 
				\\
				& \quad  
				+ \int_{0}^{\ell} \left( \left( \frac{m}{2} + \frac{\varepsilon_1 \gamma_0}{2} \right) \zeta_1 + \zeta_2 a(x) \int_{\tau_1}^{\tau_2} \, |\gamma(s)| \, ds \right) y_t^2 \, dx 
				\\ 
				& \quad 
				+ \zeta_1 \left[ -EI +  \frac{\ell^4}{4\varepsilon_1} \left( \frac{\gamma_0}{4} + \lVert a \rVert_{L^\infty(0,\ell)} \int_{\tau_1}^{\tau_2} \big| \gamma(s) \big| \, ds \right) \right] \int_0^\ell y_{xx}^2 \, dx 
				\\
				& \quad 
				- e^{-\tau_2} \zeta_2 \int_{0}^{\ell} a(x) \int_{\tau_1}^{\tau_2} \, s|\gamma(s)| \int_{0}^{1} y_{t}^2(x,t-s\varrho) \, d\varrho \, ds \, dx 
				\\
				& \quad 
				+ \left( \frac{\varepsilon_1}{2} \zeta_1 - e^{-\tau_2} \zeta_2 \right) \int_{0}^{\ell} a(x) \int_{\tau_1}^{\tau_2} \big| \gamma(s) \big| \, y_t^2(x,t-s) \, ds \, dx
				\\
				& \leqslant 
				\int_{0}^{\ell} \bigg[ -\left( \gamma_0 - \lVert a \rVert_{L^\infty(0,\ell)} \int_{\tau_1}^{\tau_2} |\gamma(s)| \, ds \right) 
				\\
				& \quad + \zeta_1 \left( \frac{m}{2} + \frac{\varepsilon_1 \gamma_0}{2} \right) + \lVert a \rVert_{L^\infty(0,\ell)} \, \zeta_2 \int_{\tau_1}^{\tau_2} \, |\gamma(s)| \, ds \bigg] y_t^2 \, dx
				\\ 
				& \quad 
				+ \zeta_1 \left[ -EI +  \frac{\ell^4}{4\varepsilon_1} \left( \frac{\gamma_0}{4} + \lVert a \rVert_{L^\infty(0,\ell)} \int_{\tau_1}^{\tau_2} \big| \gamma(s) \big| \, ds \right) \right] \int_0^\ell y_{xx}^2 \, dx 
				\\
				& \quad 
				- e^{-\tau_2} \zeta_2 \int_{0}^{\ell} a(x) \int_{\tau_1}^{\tau_2} \, s|\gamma(s)| \int_{0}^{1} y_{t}^2(x,t-s\varrho) \, d\varrho \, ds \, dx 
				\\
				& \quad 
				+ \left( \frac{\varepsilon_1}{2} \zeta_1 - e^{-\tau_2} \zeta_2 \right) \int_{0}^{\ell} a(x) \int_{\tau_1}^{\tau_2} \big| \gamma(s) \big| \, y_t^2(x,t-s) \, ds \, dx.
			\end{aligned}
		\end{gather} 
		Setting
		\begin{gather}
			\varepsilon_1 
			= 
			\frac{\ell^4}{2 EI} \left( \frac{\gamma_0}{4} + \lVert a \rVert_{L^\infty(0,\ell)} \int_{\tau_1}^{\tau_2} \big| \gamma(s) \big| \, ds \right),
		\end{gather}
		and choosing $\zeta_1, \zeta_2 > 0$ sufficiently small such that
		\begin{gather}
			\left\{
			\begin{aligned}
				& \frac{\varepsilon_1}{2} \zeta_1 - e^{-\tau_2} \zeta_2 \leqslant 0 
				\\
				& \left( \frac{m}{2} + \frac{\varepsilon_1 \gamma_0}{2} \right) \zeta_1 + \lVert a \rVert_{L^\infty(0,\ell)} \int_{\tau_1}^{\tau_2} \, |\gamma(s)| \, ds \, \zeta_2 < \gamma_0 - \lVert a \rVert_{L^\infty(0,\ell)} \int_{\tau_1}^{\tau_2} |\gamma(s)| \, ds
			\end{aligned}
			\right.
		\end{gather}
		the estimate \eqref{proof-t2-1} reduces to:
		\begin{gather}
			\begin{aligned}
				& \mathscr{L}'(t) 
				\leqslant 
				\int_{0}^{\ell} \bigg[ -\left( \gamma_0 - \lVert a \rVert_{L^\infty(0,\ell)} \int_{\tau_1}^{\tau_2} |\gamma(s)| \, ds \right) 
				\\
				& \qquad \qquad + \zeta_1 \left( \frac{m}{2} + \frac{\varepsilon_1 \gamma_0}{2} \right) + \lVert a \rVert_{L^\infty(0,\ell)} \, \zeta_2 \int_{\tau_1}^{\tau_2} \, |\gamma(s)| \, ds \bigg] y_t^2 \, dx
				\\
				& \quad 
				- \zeta_1 \frac{EI}{2} \int_{0}^{\ell} y_{xx}^2 \, dx  - e^{-\tau_2} \zeta_2 \int_{0}^{\ell} a(x) \int_{\tau_1}^{\tau_2} \, s|\gamma(s)| \int_{0}^{1} y_{t}^2(x,t-s\varrho) \, d\varrho \, ds \, dx.
			\end{aligned}
		\end{gather} 
		This implies that
		\begin{equation} \label{proof-t2-2}
			\mathscr{L}'(t) \leqslant - C E(t),
		\end{equation}
		where $C > 0$ is a constant defined by: 
		\begin{gather}
			\begin{aligned}
				& C
				:= \min \bigg\{
				\zeta_1, \, 2 e^{-\tau_2} \zeta_2, \, \mu_3 
				\bigg\},
			\end{aligned}
			\\
			\begin{aligned}
				& \mu_3:= 
				2 m^{-1} \left[ \left( \gamma_0 - \lVert a \rVert_{L^\infty(0,\ell)} \int_{\tau_1}^{\tau_2} |\gamma(s)| \, ds \right)  \right.
				\\
				& \qquad \left. - \zeta_1 \left( \frac{m}{2} + \frac{\varepsilon_1 \gamma_0}{2} \right) - \lVert a \rVert_{L^\infty(0,\ell)} \, \zeta_2 \int_{\tau_1}^{\tau_2} \, |\gamma(s)| \, ds \right].
			\end{aligned}
		\end{gather}
		Combining inequality \eqref{proof-t2-2} with the equivalence relation from Lemma \ref{l7} yields
		\begin{gather}
			\mathscr{L}'(t) \leqslant - \frac{C}{\mu_1} \mathscr{L}(t).
		\end{gather}
		Then, applying Gronwall's lemma to this differential inequality and invoking Lemma \ref{l7} once more directly establishes the desired decay estimate \eqref{ineq:expo-decay}.
	\end{proof}
	
	\bibliographystyle{plain}
	\bibliography{biblio}
	
\end{document}